\documentclass[a4paper,12pt,oneside]{article}
\usepackage{amsmath,amsfonts,amssymb,amsmath,latexsym,amsthm,enumerate}

\newtheorem{thm}{Theorem}
\newtheorem{cor}[thm]{Corollary}
\theoremstyle{definition}

\theoremstyle{plain}

\begin{document}
\title {Higher-order Cauchy of the second kind and poly-Cauchy of the second kind mixed type polynomials}
\author{by \\Dae San Kim and Taekyun Kim}\date{}\maketitle

\begin{abstract}
\noindent In this paper, we investigate some properties of higher-order Cauchy of the second kind and poly-Cauchy of the second mixed type polynomials with umbral calculus viewpoint. From our investigation, we derive many interesting identities of higher-order Cauchy of the second kind and poly-Cauchy of the second kind mixed type polynomials.
\end{abstract}

\section{Introduction}

For $\lambda\in\mathbf{C}$ with $\lambda\neq 1$, the Frobenius-Euler polynomials of order $\alpha$ ($\alpha\in\mathbf{N}\cup\{0\}$) are defined by the generating function to be
\begin{equation}\label{eq:1}
\left(\frac{1-\lambda}{e^{t}-\lambda}\right)^{\alpha}e^{xt}=\sum_{n=0}^{\infty}H_{n}^{(\alpha)}\left(x\vert\lambda\right) \frac{t^{n}}{n!},\,\,\,\,(\text{see}\,\,\lbrack 2,4,5,6,13\rbrack).
\end{equation}
As is well known, the Bernoulli polynomials of order $\alpha\in\mathbf{N}\cup\{0\}$ are also defined by the generating function to be
\begin{equation}\label{eq:2}
\left(\frac{t}{e^{t}-1}\right)^{\alpha}e^{xt}=\sum_{n=0}^{\infty}B_{n}^{(\alpha)}(x)\frac{t^{n}}{n!},\,\,\,\,(\text{see}\,\,\lbrack 1,2,4,7,8\rbrack).
\end{equation}
When $x=0$, $B_{n}^{(\alpha)}=B_{n}^{(\alpha)}(0)$ are called the Bernoulli numbers of order $\alpha$.\\
The Stirling number of the first kind is defined by
\begin{equation}\label{eq:3}
(x)_{n}=x(x-1)\cdots(x-n+1)=\sum_{l=0}^{n}S_{1}(n,l)x^{l},\,\,\,\,(n\in\mathbf{N}\cup\{0\}).
\end{equation}
From (\ref{eq:3}), we note that
\begin{equation}\label{eq:4}
\left(\log{(1+x)}\right)^{m}=m!\sum_{l=m}^{\infty}S_{1}(l,m)\frac{x^{l}}{l!}=\sum_{l=0}^{\infty}\frac{m!}{(l+m)!}S_{1}(l+m,m)x^{l+m}.
\end{equation}
It is known that the Stirling number of the second kind is given by
\begin{equation}\label{eq:5}
\left(e^{x}-1\right)^{m}=m!\sum_{l=m}^{\infty}\frac{S_{2}(l+m)}{l!}x^{l},\,\,\,\,(m\in\mathbf{N}\cup\{0\}),\,\,\,\,(\text{see}\,\,\lbrack 14,15\rbrack).
\end{equation}
The poly-logarithm factorial function is defined by
\begin{equation}\label{eq:6}
Li{f_{k}}(x)=\sum_{m=0}^{\infty}\frac{x^{m}}{m!(m+1)^{k}},\,\,\,\,(k\in\mathbf{Z}),\,\,\,\,(\text{see}\,\,\lbrack 9,10,11\rbrack).
\end{equation}
The poly-Cauchy polynomials of the second kind is given by
\begin{equation}\label{eq:7}
Li{f_{k}}\left(-\log{(1+x)}\right)(1+t)^{x}=\sum_{n=0}^{\infty}\tilde{C}_{n}^{(k)}(x)\frac{t^{n}}{n!},\,\,\,\,(\text{see}\,\,\lbrack 10-12\rbrack).
\end{equation}
and the Cauchy numbers of the second kind with order $r$ ($r\in\mathbf{N}\cup\{0\}$) are defined by the generating function to be
\begin{equation}\label{eq:8}
\left(\frac{t}{(1+t)\log{(1+t)}}\right)^{r}=\sum_{n=0}^{\infty}\mathbb{C}_{n}^{(r)}\frac{t^{n}}{n!},\,\,\,\,(\text{see}\,\,\lbrack 10-12\rbrack).
\end{equation}
Now, we consider the polynomials $\tilde{A}_{n}^{(r,k)}(x)$ whose generating function is defined by
\begin{equation}\label{eq:9}
\left(\frac{t}{(1+t)\log{(1+t)}}\right)^{r}Li{f_{k}}\left(-\log{(1+t)}\right)(1+t)^{x}=\sum_{n=0}^{\infty}\tilde{A}_{n}^{(r,k)}\frac{t^{n}}{n!},
\end{equation}
where $r\in\mathbf{N}\cup\{0\}$ and $k\in\mathbf{Z}$.\\
$\tilde{A}_{n}^{(r,k)}(x)$ are called higher-order Cauchy of the second kind and poly-Cauchy of the second kind mixed type polynomials. When $x=0$, $\tilde{A}_{n}^{(r,k)}=\tilde{A}_{n}^{(r,k)}(0)$  are called the higher-order Cauchy of the second kind and poly-Cauchy of the second kind mixed type numbers.\\
Let $\mathcal{F}$ be the set of all formal power series in the variable $t$ over $\mathbf{C}$ as follows:
\begin{equation}\label{eq:10}
\mathcal{F}=\left\{f(t)=\sum_{k=0}^{\infty}a_{k}\frac{t^{k}}{k!}\Bigg\vert a_{k}\in\mathbf{C}\right\}.
\end{equation}
Let $\mathbb{P}=\mathbf{C}\lbrack x\rbrack$ and let $\mathbb{P}^{*}$ be the vector space of all linear functionals on $\mathbb{P}$. $\langle L\vert p(x)\rangle$ denotes the action of the linear functional $L$ on the polynomial $p(x)$, and the vector space operations on $\mathbb{P}^{*}$ are defined by $\langle L+M\vert p(x)\rangle=\langle L\vert p(x)\rangle +\langle M\vert p(x)\rangle$, $\langle cL\vert p(x)\rangle=c\langle L\vert p(x)\rangle$, where $c$ is a complex constant. For $f(t)\in\mathcal{F}$, let $\langle f(t)\vert x^{n}\rangle =a_{n}$. Then, by (\ref{eq:10}), we see that
\begin{equation}\label{eq:11}
\langle t^{k}\vert x^{n}\rangle=n!\delta_{n,k},\,\,\,\,(\text{see}\,\,\lbrack 14,15\rbrack),
\end{equation}
where $\delta_{n,k}$ is the Kronecker's symbol.\\
Let us assume that $f_{L}(t)=\sum_{k=0}^{\infty}\langle L\vert x^{k}\rangle\frac{t^{k}}{k!}$. By (\ref{eq:11}), we see that $\langle f_{L}(t)\vert x^{n}\rangle=\langle L\vert x^{n}\rangle$. That is, $f_{L}(t)=L$. Additionally, the map $L\longmapsto f_{L}(t)$ is a vector space isomorphism from $\mathbb{P}^{*}$ onto $\mathcal{F}$. Henceforth, $\mathcal{F}$ denotes both the algebra of the formal power series in $t$ and the vector space of all linear functionals on $\mathbb{P}$, and so an element $f(t)$ of $\mathcal{F}$ will be thought as a formal power series and a linear functional. $\mathcal{F}$ is called the umbral algebra. The umbral calculus is the study of umbral algebra.\\
The order of the power series $f(t)(\neq 0)$ is the smallest integer for which $a_{k}$ does not vanish. The order of $f(t)$ is denoted by $O\left(f(t)\right)$. If $O\left(f(t)\right)=0$, then $f(t)$ is called an invertible series. If $O\left(f(t)\right)=1$, then $f(t)$ is said to be a delta series.\\
For $f(t)\in\mathcal{F}$ and $p(x)\in\mathbb{P}$, we have
\begin{equation}\label{eq:12}
f(t)=\sum_{k=0}^{\infty}\langle f(t)\vert x^{k}\rangle\frac{t^{k}}{k!},\,\,\,\, p(x)=\sum_{k=0}^{\infty}\langle t^{k}\vert p(x)\rangle\frac{x^{k}}{k!},\,\,\,\,(\text{see}\,\,\lbrack 14\rbrack).
\end{equation}
Thus, by (\ref{eq:12}), we get
\begin{equation}\label{eq:13}
p^{(k)}(x)=\frac{d^{k}p(x)}{dx^{k}}=\sum_{l=k}^{\infty}\frac{\langle t^{l}\vert p(x)\rangle}{(l-k)!}x^{l-k},\,\,\,\, p^{(k)}(0)=\langle t^{k}\vert p(x)\rangle=\langle 1\vert p^{(k)}(x)\rangle.
\end{equation}
From (\ref{eq:13}), we note that
\begin{equation}\label{eq:14}
t^{k}p(x)=p^{(k)}(x),\,\,\,\,e^{yt}p(x)=p(x+y),\,\,\,\,\langle e^{yt}\vert p(x)\rangle=p(y).
\end{equation}
For $O\left(f(t)\right)=1$, $O\left(g(t)\right)=0$, there exists a unique sequence $s_{n}(x)$ of polynomials such that $\langle g(t)f(t)^{k}\vert s_{n}(x)\rangle=n!\delta_{n,k}$, for $n, k\geq 0$. The sequence $s_{n}(x)$ is called the Sheffer sequence for $\left(g(t), f(t)\right)$ which is denoted by $s_{n}(x)\sim\left(g(t), f(t)\right)$.\\
Let $p(x)\in\mathbb{P}$ and $f(t)\in\mathcal{F}$. Then we see that
\begin{equation}\label{eq:15}
\langle f(t)\vert xp(x)\rangle=\langle \partial_{t}f(t)\vert p(x)\rangle=\langle f'(t)\vert p(x)\rangle,\,\,\,\,(\text{see}\,\,\lbrack 14\rbrack).
\end{equation}
For $s_{n}(x)\sim\left(g(t),f(t)\right)$, we have the following equations:
\begin{equation}\label{eq:16}
s_{n}(x)=\sum_{j=0}^{n}\frac{1}{j!}\left\langle g\left(\bar{f}(t)\right)^{-1}\bar{f}(t)^{j}\Big\vert x^{n}\right\rangle x^{j},
\end{equation}
where $\bar{f}(t)$ is the compositional inverse of $f(t)$ with $\bar{f}\left(f(t)\right)=t$,
\begin{equation}\label{eq:17}
\frac{1}{g\left(\bar{f}(t)\right)}e^{x\bar{f}(t)}=\sum_{n=0}^{\infty}s_{n}(x)\frac{t^{n}}{n!},\,\,\,\,\text{for all}\,\, x\in\mathbf{C},
\end{equation}
\begin{equation}\label{eq:18}
s_{n}(x+y)=\sum_{k=0}^{n}\binom{n}{k}s_{k}(x)p_{n-k}(y),\,\,\,\,\text{where}\,\,p_{n-k}(y)=g(t)s_{n-k}(y),
\end{equation}
and
\begin{equation}\label{eq:19}
s_{n+1}(x)=\left(x-\frac{g'(t)}{g(t)}\right)\frac{1}{f'(t)}s_{n}(x),\,\,\,\, f(t)s_{n}(x)=ns_{n-1}(x),\,\,\,\,(\text{see}\,\,\lbrack 3,5,9,14\rbrack).
\end{equation}
For $s_{n}(x)\sim\left(g(t),f(t)\right)$, $r_{n}(x)\sim\left(h(t),l(t)\right)$, we have
\begin{equation}\label{eq:20}
s_{n}(x)=\sum_{m=0}^{n}C_{n,m}r_{m}(x),
\end{equation}
where
\begin{equation}\label{eq:21}
C_{n,m}=\frac{1}{m!}\left\langle\frac{h\left(\bar{f}(t)\right)}{g\left(\bar{f}(t)\right)}l\left(\bar{f}(t)\right)^{m}\Bigg\vert x^{n}\right\rangle,\,\,\,\,(\text{see}\,\,\lbrack 14\rbrack).
\end{equation}
In this paper, we consider higher-order Cauchy of the second kind and poly-Cauchy of the second kind mixed type polynomials and we investigate some properties of those polynomials with umbral calculus viewpoint. From our investigation, we can derive many interesting identities related to higher-order Cauchy of the second kind and poly-Cauchy of the second kind mixed type polynomials.

\section{Higer-order Cauchy of the second kind and poly-Cauchy of the second kind mixed type polynomials}

From (\ref{eq:3}) and (\ref{eq:17}), we note that $\tilde{A}_{n}^{(r,k)}(x)$ is the Sheffer sequence for the pair $\left(\left(\frac{te^{t}}{e^{t}-1}\right)^{r}\frac{1}{Li{f_{k}}(-t)}, e^{t}-1\right)$. That is,
\begin{equation}\label{eq:22}
\tilde{A}_{n}^{(r,k)}(x)\sim\left(\left(\frac{te^{t}}{e^{t}-1}\right)^{r}\frac{1}{Li{f_{k}}(-t)}, e^{t}-1\right).
\end{equation}
Komatsu considered the number $\tilde{A}_{n}^{(r,k)}$, which was denoted by $\tilde{T}_{r+1}^{(k)}(n)$ (see $\lbrack 10-12\rbrack$).\\
By (\ref{eq:22}), we easily see that
\begin{equation}\label{eq:23}
\left(\frac{te^{t}}{e^{t}-1}\right)^{r}\frac{1}{Li{f_{k}}(-t)}\tilde{A}_{n}^{(r,k)}(x)\sim\left(1, e^{t}-1\right),
\end{equation}
and we see that $(x)_{n}\sim\left(1,e^{t}-1\right)$.\\
From the uniqueness of Sheffer sequence, we note that
\begin{equation}\label{eq:24}
\left(\frac{te^{t}}{e^{t}-1}\right)^{r}\frac{1}{Li{f_{k}}(-t)}A_{n}^{(r,k)}(x)=(x)_{n}=\sum_{m=0}^{n}S_{1}(n,m)x^{m}\sim\left(1,e^{t}-1\right).
\end{equation}
By (\ref{eq:24}), we get
\begin{align}\label{eq:25}
A_{n}^{(r,k)}(x)&=\left(\frac{e^{t}-1}{te^{t}}\right)^{r}Li{f_{k}}(-t)(x)_{n}=\sum_{m=0}^{n}S_{1}(n,m)\left(\frac{e^{t}-1}{te^{t}}\right)^{r}Li{f_{k}}(-t)x^{m}\\
&=\sum_{m=0}^{n}S_{1}(n,m)\left(\frac{e^{t}-1}{te^{t}}\right)^{r}\sum_{l=0}^{\infty}\frac{(-t)^{l}}{l!(l+1)^{k}}x^{m}\nonumber\\
&=\sum_{m=0}^{n}S_{1}(n,m)\sum_{l=0}^{m}\frac{(-1)^{l}(m)_{l}}{l!(l+1)^{k}}\left(\frac{e^{-t}-1}{-t}\right)^{r}x^{m-l}\nonumber\\
&=\sum_{m=0}^{n}S_{1}(n,m)\sum_{l=0}^{m}\frac{(-1)^{l}(m)_{l}}{l!(l+1)^{k}}\nonumber\\
&\quad\times\sum_{a=0}^{m-l}\frac{r!}{(a+r)!}S_{2}(a+r,r)(-1)^{a}(m-l)_{a}x^{m-l-a}\nonumber
\end{align}
\begin{align*}
&=\sum_{m=0}^{n}\sum_{l=0}^{m}\sum_{a=0}^{m-l}(-1)^{m}\frac{\binom{m}{l}\binom{m-l}{a}}{\binom{a+r}{r}(l+1)^{k}}S_{1}(n,m)S_{2}(a+r,r)(-x)^{m-l-a}\\
&=\sum_{j=0}^{n}\left\{\sum_{m=j}^{n}\sum_{l=0}^{m-j}(-1)^{m}\frac{\binom{m}{l}\binom{m-l}{j}}{\binom{m-l-j+r}{r}(l+1)^{k}}S_{1}(n,m)S_{2}(m-l-j+r,r)\right\}(-x)^{j}.
\end{align*}
Therefore, by (\ref{eq:25}), we obtain the following theorem.

\begin{thm}\label{eq:thm1}
For $n, r\geq 0$, $k\in\mathbf{Z}$, we have
\begin{align*}
\tilde{A}_{n}^{(r,k)}(x)&=\sum_{j=0}^{n}\Bigg\{\sum_{m=j}^{n}\sum_{l=0}^{m-j}(-1)^{m}\frac{\binom{m}{l}\binom{m-l}{j}}{\binom{m-l-j+r}{r}(l+1)^{k}}\\
&\quad\times S_{1}(n,m)S_{2}(m-l-j+r,r)\Bigg\}(-x)^{j}.
\end{align*}
\end{thm}

\noindent From (\ref{eq:16}) and (\ref{eq:22}), we note that
\begin{align}\label{eq:26}
&\tilde{A}_{n}^{(r,k)}(x)\\
&=\sum_{j=0}^{n}\frac{1}{j!}\left\langle\left(\frac{t}{(1+t)\log{(1+t)}}\right)^{r}Li{f_{k}}\left(-\log{(1+t)}\right)\left(\log{(1+t)}\right)^{j} \Bigg\vert x^{n}\right\rangle x^{j}\nonumber\\
&=\sum_{j=0}^{n}\frac{1}{j!}\sum_{l=0}^{n-j}\frac{j!}{(l+j)!}S_{1}(l+j,j)(n)_{l+j}\nonumber\\
&\quad\times\left\langle\left(\frac{t}{(1+t)\log{(1+t)}}\right)^{r}Li{f_{k}}\left(-\log{(1+t)}\right)\Bigg\vert x^{n-l-j}\right\rangle x^{j},\nonumber
\end{align}
and
\begin{align}\label{eq:27}
&\left\langle\left(\frac{t}{(1+t)\log{(1+t)}}\right)^{r}Li{f_{k}}\left(-\log{(1+t)}\right)\Bigg\vert x^{n-l-j}\right\rangle\\
&=\sum_{a=0}^{\infty}\frac{\tilde{A}_{a}^{(r,k)}}{a!}\left\langle t^{a}\big\vert x^{n-l-j}\right\rangle\nonumber\\
&=\tilde{A}_{n-l-j}^{(r,k)}.\nonumber
\end{align}
By (\ref{eq:26}) and (\ref{eq:27}), we get
\begin{align}\label{eq:28}
\tilde{A}_{n}^{(r,k)}(x)&=\sum_{j=0}^{n}\left\{\sum_{l=0}^{n-j}\binom{n}{l+j}
S_{1}(l+j,j)\tilde{A}_{n-l-j}^{(r,k)}\right\}x^{j}\\
&=\sum_{j=0}^{n}\left\{\sum_{l=0}^{n-j}\binom{n}{l}S_{1}(n-l,j)\tilde{A}_{l}^{(r,k)}\right\}x^{j}.\nonumber
\end{align}
Therefore, by (\ref{eq:28}), we obtain the following theorem.

\begin{thm}\label{eq:thm2}
For $n, r\geq 0$ and $k\in\mathbf{Z}$, we have
\begin{equation*}
\tilde{A}_{n}^{(r,k)}(x)=\sum_{j=0}^{n}\left\{\sum_{l=0}^{n-j}\binom{n}{l}S_{1}(n-l,j)\tilde{A}_{l}^{(r,k)}\right\}x^{j}.
\end{equation*}
\end{thm}

\noindent It is known that
\begin{equation}\label{eq:29}
\left(\frac{t}{\log{(1+t)}}\right)^{n}(1+t)^{x-1}=\sum_{k=0}^{\infty}B_{k}^{(k-n+1)}(x)\frac{t^{k}}{k!}.
\end{equation}
In particular, for $x=1-r$, $n=r$, we have
\begin{equation}\label{eq:30}
\left(\frac{t}{(1+t)\log{(1+t)}}\right)^{r}=\sum_{k=0}^{\infty}B_{k}^{(k-r+1)}(1-r)\frac{t^{k}}{k!}.
\end{equation}
By (\ref{eq:26}) and (\ref{eq:30}), we get
\begin{align}\label{eq:31}
\tilde{A}_{n}^{(r,k)}(x)&=\sum_{j=0}^{n}\sum_{l=0}^{n-j}\binom{n}{l+j}S_{1}(l+j,j)\sum_{a=0}^{\infty}\frac{B_{a}^{(a-r+1)}(1-r)}{a!}\\
&\quad\times\left\langle Li{f_{k}}\left(-\log{(1+t)}\right)\big\vert t^{a}x^{n-l-j}\right\rangle x^{j}\nonumber\\
&=\sum_{j=0}^{n}\sum_{l=0}^{n-j}\binom{n}{l+j}S_{1}(l+j,j)\sum_{a=0}^{n-l-j}B_{a}^{(a-r+1)}(1-r)\frac{1}{a!}(n-l-j)_{a}\nonumber\\
&\quad\times \left\langle Li{f_{k}}\left(-\log{(1+t)}\right)\big\vert x^{n-l-j-a}\right\rangle x^{j}\nonumber\\
&=\sum_{j=0}^{n}\Bigg\{\sum_{l=0}^{n-j}\sum_{a=0}^{n-j-l}\binom{n}{l+j}\binom{n-j-l}{a}S_{1}(l+j,j)B_{a}^{(a-r+1)}(1-r)\nonumber\\
&\quad\times\tilde{C}_{n-j-l-a}^{(k)}\Bigg\}x^{j},\nonumber
\end{align}
where $\tilde{C}_{n}^{(k)}$ are the poly-Cauchy numbers of the second kind.\\
Therefore, by (\ref{eq:31}), we obtain the following theorem.

\begin{thm}\label{eq:thm3}
For $n,r\geq 0$ and $k\in\mathbf{Z}$, we have
\begin{align*}
\tilde{A}_{n}^{(r,k)}(x)&=\sum_{j=0}^{n}\Bigg\{\sum_{l=0}^{n-j}\sum_{a=0}^{n-j-l}\binom{n}{l+j}\binom{n-j-l}{a}S_{1}(l+j,j)B_{a}^{(a-r+1)}(1-r)\\
&\quad\times\tilde{C}_{n-j-l-a}^{(k)}\Bigg\}x^{j}.
\end{align*}
\end{thm}

\noindent By (\ref{eq:29}), we easily see that
\begin{equation}\label{eq:32}
\frac{t}{(1+t)\log{(1+t)}}=\sum_{n=0}^{\infty}B_{n}^{(n)}\frac{t^{n}}{n!}.
\end{equation}
Thus, by (\ref{eq:26}) and (\ref{eq:32}), we get
\begin{align}\label{eq:33}
\tilde{A}_{n}^{(r,k)}&=\sum_{j=0}^{n}\sum_{l=0}^{n-j}\binom{n}{l+j}S_{1}(l+j,j)\\
&\quad\times\left\langle Li{f_{k}}\left(-\log{(1+t)}\right)\Big\vert \left(\frac{t}{(1+t)\log{(1+t)}}\right)^{r}x^{n-l-j}\right\rangle x^{j}\nonumber\\
&=\sum_{j=0}^{n}\Bigg\{\sum_{l=0}^{n-j}\sum_{a=0}^{n-j-l}\sum_{a_{1}+\cdots+a_{r}=a}\binom{n}{l+j}\binom{n-j-l}{a}\binom{a}{a_{1},\cdots,a_{r}}S_{1}(l+j,j)\nonumber\\
&\quad\times \left(\prod_{i=1}^{r}B_{a_{i}}^{(a_{i})}\right)\tilde{C}_{n-j-l-a}^{(k)}\Bigg\}x^{j}.\nonumber
\end{align}
Therefore, by (\ref{eq:33}), we obtain the following corollary.

\begin{cor}\label{eq:cor4}
For $n\geq 0$, $r\in\mathbf{N}$ and $k\in\mathbf{Z}$, we have
\begin{align*}
\tilde{A}_{n}^{(r,k)}(x)&=\sum_{j=0}^{n}\Bigg\{\sum_{l=0}^{n-j}\sum_{a=0}^{n-j-l}\sum_{a_{1}+\cdots+a_{r}=a}\binom{n}{l+j}\binom{n-j-l}{a}\binom{a}{a_{1},\cdots,a_{r}}\nonumber\\
&\quad\times S_{1}(l+j,j)\left(\prod_{i=1}^{r}B_{a_{i}}^{(a_{i})}\right)\tilde{C}_{n-j-l-a}^{(k)}\Bigg\}x^{j}.\nonumber
\end{align*}
\end{cor}

\noindent From (\ref{eq:18}) and (\ref{eq:19}), we can derive
\begin{equation}\label{eq:34}
\tilde{A}_{n}^{(r,k)}(x+y)=\sum_{j=0}^{n}\binom{n}{j}\tilde{A}_{j}^{(r,k)}(x)(y)_{n-j},\,\,\,\, \left(e^{t}-1\right)\tilde{A}_{n}^{(r,k)}(x)=n\tilde{A}_{n-1}^{(r,k)}(x).
\end{equation}
By (\ref{eq:19}) and (\ref{eq:22}), we get
\begin{align}\label{eq:35}
\tilde{A}_{n+1}^{(r,k)}(x)&=x\tilde{A}_{n}^{(r,k)}(x-1)-r\sum_{m=0}^{n}\sum_{l=0}^{m}\sum_{a=0}^{m-l}\frac{(-1)^{m-a}\binom{m}{l}\binom{m-l}{a}}{(a+2)(a+1)(l+1)^{k}}S_{1}(n,m)\\
&\quad\times B_{m-l-a}^{(1-r)}(2-x)-\sum_{m=0}^{n}\sum_{a=0}^{m}\frac{(-1)^{m}\binom{m}{a}}{(a+2)^{k}}S_{1}(n,m)B_{m-a}^{(-r)}(1-x).\nonumber
\end{align}
From (\ref{eq:11}), we note that
\begin{align}\label{eq:36}
\tilde{A}_{n}^{(r,k)}(y)&=\left\langle\sum_{m=0}^{\infty}\tilde{A}_{m}^{(r,k)}(y)\frac{t^{m}}{m!}\Bigg\vert x^{n}\right\rangle\\
&=\left\langle\left(\frac{t}{(1+t)\log{(1+t)}}\right)^{r}Li{f_{k}}\left(-\log{(1+t)}\right)(1+t)^{y}\Bigg\vert x^{n}\right\rangle\nonumber\\
&=\left\langle\left(\frac{t}{(1+t)\log{(1+t)}}\right)^{r}Li{f_{k}}\left(-\log{(1+t)}\right)(1+t)^{y}\Bigg\vert xx^{n-1}\right\rangle.\nonumber
\end{align}
By (\ref{eq:15}) and (\ref{eq:36}), we get
\begin{align}\label{eq:37}
\tilde{A}_{n}^{(r,k)}(y)&=\left\langle\partial_{t}\left(\left(\frac{t}{(1+t)\log{(1+t)}}\right)^{r}Li{f_{k}}\left(-\log{(1+t)}\right)(1+t)^{y}\right)\Bigg\vert x^{n-1}\right\rangle\\
&=\left\langle\left(\partial_{t}\left(\frac{t}{(1+t)\log{(1+t)}}\right)^{r}\right)Li{f_{k}}\left(-\log{(1+t)}\right)(1+t)^{y}\Bigg\vert x^{n-1}\right\rangle\nonumber\\
&\quad +\left\langle\left(\frac{t}{(1+t)\log{(1+t)}}\right)^{r}\left(\partial_{t}\left(Li{f_{k}}\left(-\log{(1+t)}\right)\right)\right)(1+t)^{y}\Bigg\vert x^{n-1}\right\rangle\nonumber\\
&\quad +\left\langle\left(\frac{t}{(1+t)\log{(1+t)}}\right)^{r}Li{f_{k}}\left(-\log{(1+t)}\right)\left(\partial_{t}(1+t)^{y}\right)\Bigg\vert x^{n-1}\right\rangle\nonumber\\
&=y\tilde{A}_{n-1}^{(r,k)}(y-1)\nonumber\\
&\quad +\left\langle\left(\partial_{t}\left(\frac{t}{(1+t)\log{(1+t)}}\right)^{r}\right)Li{f_{k}}\left(-\log{(1+t)}\right)(1+t)^{y}\Bigg\vert x^{n-1}\right\rangle\nonumber\\
&\quad +\left\langle\left(\frac{t}{(1+t)\log{(1+t)}}\right)^{r}\left(\partial_{t}\left(Li{f_{k}}\left(-\log{(1+t)}\right)\right)\right) (1+t)^{y}\Bigg\vert x^{n-1}\right\rangle.\nonumber
\end{align}
Now, we observe that
\begin{align}\label{eq:38}
\frac{\log{(1+t)}-t}{t^{2}}x^{l}&=\sum_{a=0}^{l}\frac{(-1)^{a-1}}{a+2}t^{a}x^{l}=\sum_{a=0}^{l}\frac{(-1)^{a-1}}{a+2}(l)_{a}x^{l-a}\\
&=\sum_{a=0}^{l}\frac{(-1)^{l-a-1}}{l-a+2}\binom{l}{a}(l-a)!x^{a}.\nonumber
\end{align}
By (\ref{eq:38}), we get
\begin{align}\label{eq:39}
&\left\langle\left(\partial_{t}\left(\frac{t}{(1+t)\log{(1+t)}}\right)^{r}\right)Li{f_{k}}\left(-\log{(1+t)}\right)(1+t)^{y}\Bigg\vert x^{n-1}\right\rangle\\
&=r\left\langle\left(\frac{t}{(1+t)\log{(1+t)}}\right)^{r+1}\frac{\log{(1+t)}-t}{t^{2}}Li{f_{k}}\left(-\log{(1+t)}\right)(1+t)^{y}\Bigg\vert x^{n-1}\right\rangle\nonumber\\
&=r\sum_{a=0}^{n-1}\frac{(-1)^{n-a}(n-a-1)!}{n-a+1}\binom{n-1}{a}\nonumber\\
&\quad\times\left\langle\left(\frac{t}{(1+t)\log{(1+t)}}\right)^{r+1}Li{f_{k}}\left(-\log{(1+t)}\right)(1+t)^{y}\Bigg\vert x^{a}\right\rangle\nonumber\\
&=r\sum_{a=0}^{n-1}\frac{(-1)^{n-a}(n-a-1)!}{n-a+1}\binom{n-1}{a}\tilde{A}_{a}^{(r+1,k)}(y)\nonumber\\
&=r\sum_{a=0}^{n-1}\frac{(-1)^{a+1}a!}{a+2}\binom{n-1}{a}\tilde{A}_{n-1-a}^{(r+1,k)}(y).\nonumber
\end{align}
It is not difficult to show that
\begin{equation}\label{eq:40}
\left(Li{f_{k}}\left(-\log{(1+t)}\right)\right)'=\frac{Li{f_{k-1}}\left(-\log{(1+t)}\right)-Li{f_{k}}\left(-\log{(1+t)}\right)}{(1+t)\log{(1+t)}}.
\end{equation}
By (\ref{eq:40}), we get
\begin{align}\label{eq:41}
&\left\langle\left(\frac{t}{(1+t)\log{(1+t)}}\right)^{r}\left(\partial_{t}Li{f_{k-1}}\left(-\log{(1+t)}\right)\right)(1+t)^{y}\Bigg\vert x^{n-1}\right\rangle\\
&=\Bigg\langle\left(\frac{t}{(1+t)\log{(1+t)}}\right)^{r}\frac{Li{f_{k-1}}\left(-\log{(1+t)}\right)-Li{f_{k}}\left(-\log{(1+t)}\right)}{(1+t)\log{(1+t)}}\nonumber\\
&\quad\times(1+t)^{y}\Bigg\vert \frac{1}{n}tx^{n}\Bigg\rangle\nonumber
\end{align}
\begin{align*}
&=\frac{1}{n}\left(\tilde{A}_{n}^{(r+1,k-1)}(y)-\tilde{A}_{n}^{(r+1,k)}(y)\right).
\end{align*}
Therefore, by (\ref{eq:37}), (\ref{eq:39}) and (\ref{eq:41}), we obtain the following theorem.

\begin{thm}\label{eq:thm5}
For $n\geq 1$, $r,k\in\mathbf{Z}$ with $r\geq 1$, we have
\begin{align*}
\tilde{A}_{n}^{(r,k)}(x)&=x\tilde{A}_{n-1}^{(r,k)}(x-1)+r\sum_{a=0}^{n-1}\frac{(-1)^{a+1}a!}{a+2}\binom{n-1}{a}\tilde{A}_{n-1-a}^{(r+1,k)}(x)\\
&\quad +\frac{1}{n}\left(\tilde{A}_{n}^{(r+1,k-1)}(x)-\tilde{A}_{n}^{(r+1,k)}(x)\right).
\end{align*}
\end{thm}

\noindent Now, we compute the following equation (\ref{eq:42}) in two different ways:
\begin{align}\label{eq:42}
\left\langle\left(\frac{t}{(1+t)\log{(1+t)}}\right)^{r}Li{f_{k}}\left(-\log{(1+t)}\right)\left(\log{(1+t)}\right)^{m}\Bigg\vert x^{n}\right\rangle.
\end{align}
On the one hand,
\begin{align}\label{eq:43}
&\left\langle\left(\frac{t}{(1+t)\log{(1+t)}}\right)^{r}Li{f_{k}}\left(-\log{(1+t)}\right)\left(\log{(1+t)}\right)^{m}\Bigg\vert x^{n}\right\rangle\\
&=\sum_{l=0}^{n-m}\frac{m!}{(l+m)!}S_{1}(l+m,m)(n)_{l+m}\nonumber\\
&\quad\times\left\langle\left(\frac{t}{(1+t)\log{(1+t)}}\right)^{r}Li{f_{k}}\left(-\log{(1+t)}\right)\Bigg\vert x^{n-l-m}\right\rangle\nonumber\\
&=\sum_{l=0}^{n-m}m!\binom{n}{l+m}S_{1}(l+m,m)\tilde{A}_{n-l-m}^{(r,k)}\nonumber\\
&=\sum_{l=0}^{n-m}m!\binom{n}{l}S_{1}(n-l,m)\tilde{A}_{l}^{(r,k)}.\nonumber
\end{align}
On the other hand, (\ref{eq:42}) is
\begin{align}\label{eq:44}
&\left\langle\left(\frac{t}{(1+t)\log{(1+t)}}\right)^{r}Li{f_{k}}\left(-\log{(1+t)}\right)\left(\log{(1+t)}\right)^{m}\Bigg\vert xx^{n-1}\right\rangle\\
&=\left\langle\partial_{t}\left(\left(\frac{t}{(1+t)\log{(1+t)}}\right)^{r}Li{f_{k}}\left(-\log{(1+t)}\right)\left(\log{(1+t)}\right)^{m}\right)\Bigg\vert x^{n-1}\right\rangle\nonumber
\end{align}
\begin{align*}
&=\left\langle\left(\partial_{t}\left(\frac{t}{(1+t)\log{(1+t)}}\right)^{r}\right)Li{f_{k}}\left(-\log{(1+t)}\right)\left(\log{(1+t)}\right)^{m}\Bigg\vert x^{n-1}\right\rangle\nonumber\\
&\quad +\left\langle\left(\frac{t}{(1+t)\log{(1+t)}}\right)^{r}\left(\partial_{t}Li{f_{k}}\left(-\log{(1+t)}\right)\right)\left(\log{(1+t)}\right)^{m}\Bigg\vert x^{n-1}\right\rangle\nonumber\\
&\quad +\left\langle\left(\frac{t}{(1+t)\log{(1+t)}}\right)^{r}Li{f_{k}}\left(-\log{(1+t)}\right)\left(\partial_{t}\left(\log{(1+t)}\right)^{m}\right) \Bigg\vert x^{n-1}\right\rangle
\end{align*}
Here, we observe that
\begin{align}\label{eq:45}
&\left\langle\left(\partial_{t}\left(\frac{t}{(1+t)\log{(1+t)}}\right)^{r}\right)Li{f_{k}}\left(-\log{(1+t)}\right)\left(\log{(1+t)}\right)^{m}\Bigg\vert x^{n-1}\right\rangle\\
&=r\Bigg\langle\left(\frac{t}{(1+t)\log{(1+t)}}\right)^{r+1}\frac{\log{(1+t)}-t}{t^{2}}Li{f_{k}}\left(-\log{(1+t)}\right)\nonumber\\
&\quad\times\left(\log{(1+t)}\right)^{m}\Bigg\vert x^{n-1}\Bigg\rangle\nonumber\\
&=r\sum_{l=0}^{n-m-1}\frac{m!}{(l+m)!}S_{1}(l+m,m)(n-1)_{l+m}\nonumber\\
&\quad\times\sum_{a=0}^{n-1-m-l}\frac{(-1)^{a+1}a!\binom{n-l-m-1}{a}}{a+2}\tilde{A}_{n-l-m-a-1}^{(r+1,k)}\nonumber\\
&=r\sum_{l=0}^{n-m-1}\sum_{a=0}^{n-1-l-m}\frac{(-1)^{a+1}a!m!}{a+2}\binom{n-1}{l+m}\binom{n-l-m-1}{a}\nonumber\\
&\quad\times S_{1}(l+m,m)\tilde{A}_{n-l-m-a-1}^{(r+1,k)},\nonumber
\end{align}
and
\begin{align}\label{eq:46}
&\left\langle\left(\frac{t}{(1+t)\log{(1+t)}}\right)^{r}\left(\partial_{t}Li{f_{k}}\left(-\log{(1+t)}\right)\right)\left(\log{(1+t)}\right)^{m}\Bigg\vert x^{n-1}\right\rangle\\
&=\Bigg\langle\left(\frac{t}{(1+t)\log{(1+t)}}\right)^{r}\left(\frac{Li{f_{k-1}}\left(-log{(1+t)}\right)-Li{f_{k}}\left(-log{(1+t)}\right)}{(1+t)\log{(1+t)}}\right)\nonumber\\
&\quad\times\left(\log{(1+t)}\right)^{m}\Bigg\vert x^{n-1}\Bigg\rangle\nonumber
\end{align}
\begin{align*}
&=\sum_{l=0}^{n-m-1}\frac{m!}{(l+m)!}S_{1}(l+m,m)(n-1)_{l+m}\frac{1}{n-l-m}\left(\tilde{A}_{n-l-m}^{(r+1,k-1)}-\tilde{A}_{n-l-m}^{(r+1,k)}\right)\nonumber\\
&=\sum_{l=0}^{n-m-1}\frac{m!}{n-l-m}\binom{n-1}{l+m}S_{1}(l+m,m)\left(\tilde{A}_{n-l-m}^{(r+1,k-1)}-\tilde{A}_{n-l-m}^{(r+1,k)}\right).\nonumber
\end{align*}
Finally, we easily see that
\begin{align}\label{eq:47}
&\left\langle\left(\frac{t}{(1+t)\log{(1+t)}}\right)^{r}Li{f_{k}}\left(-\log{(1+t)}\right)\left(\partial_{t}\left(\log{(1+t)}\right)^{m}\right) \Bigg\vert x^{n-1}\right\rangle\\
&=m\left\langle\left(\frac{t}{(1+t)\log{(1+t)}}\right)^{r}Li{f_{k}}\left(-\log{(1+t)}\right)(1+t)^{-1}\left(\log{(1+t)}\right)^{m-1}\Bigg\vert x^{n-1}\right\rangle\nonumber\\
&=m\sum_{l=0}^{n-m}\frac{(m-1)!}{(l+m-1)!}S_{1}(l+m-1,m-1)(n-1)_{l+m-1}\tilde{A}_{n-l-m}^{(r,k)}(-1)\nonumber\\
&=\sum_{l=0}^{n-m}m!\binom{n-1}{l+m-1}S_{1}(l+m-1,m-1)\tilde{A}_{n-l-m}^{(r,k)}(-1).\nonumber
\end{align}
Therefore, by (\ref{eq:43}), (\ref{eq:44}), (\ref{eq:45}), (\ref{eq:46}) and (\ref{eq:47}), we obtain the following theorem.

\begin{thm}\label{eq:thm6}
For $n-1\geq m\geq 1$, we have
\begin{align*}
&\sum_{l=0}^{n-m}\binom{n}{l}S_{1}(n-l,m)\tilde{A}_{l}^{(r,k)}\\
&=r\sum_{l=0}^{n-m-1}\sum_{a=0}^{n-l-m-1}\frac{(-1)^{a+1}a!}{a+2}\binom{n-1}{l+m}\binom{n-l-m-1}{a}S_{1}(l+m,m)\tilde{A}_{n-l-m-a-1}^{(r+1,k)}\nonumber\\
&\quad +\sum_{l=0}^{n-m-1}\frac{1}{n-l-m}\binom{n-1}{l+m}S_{1}(l+m,m)\left(\tilde{A}_{n-l-m}^{(r+1,k-1)}-\tilde{A}_{n-l-m}^{(r+1,k)}\right)\nonumber\\
&\quad +\sum_{l=0}^{n-m}\binom{n-1}{l+m-1}S_{1}(l+m-1,m-1)\tilde{A}_{n-l-m}^{(r,k)}(-1).\nonumber
\end{align*}
\end{thm}

\noindent For $s_{n}(x)\sim\left(g(t),f(t)\right)$, we note that
\begin{equation}\label{eq:48}
\frac{d}{dx}s_{n}(x)=\sum_{l=0}^{n-1}\binom{n}{l}\left\langle\bar{f}(t)\big\vert x^{n-l}\right\rangle s_{l}(x).
\end{equation}
From (\ref{eq:22}) and (\ref{eq:48}), we can derive the following equation (\ref{eq:49}):
\begin{align}\label{eq:49}
\frac{d}{dx}\tilde{A}_{n}^{(r,k)}(x)&=\sum_{l=0}^{n-1}\binom{n}{l}\left\langle\log{(1+t)}\big\vert x^{n-l}\right\rangle\tilde{A}_{l}^{(r,k)}(x)\\
&=\sum_{l=0}^{n-1}\binom{n}{l}\sum_{m=0}^{\infty}\frac{(-1)^{m}}{m+1}\left\langle t^{m+1}\big\vert x^{n-l}\right\rangle\tilde{A}_{l}^{(r,k)}(x)\nonumber\\
&=\sum_{l=0}^{n-1}\binom{n}{l}\frac{(-1)^{n-l-1}}{n-l}(n-l)!\tilde{A}_{l}^{(r,k)}(x)\nonumber\\
&=(-1)^{n}n!\sum_{l=0}^{n-1}\frac{(-1)^{l+1}}{(n-l)l!}\tilde{A}_{l}^{(r,k)}(x).\nonumber
\end{align}
For $\tilde{A}_{n}^{(r,k)}(x)\sim\left(\left(\frac{te^{t}}{e^{t}-1}\right)^{r}\frac{1}{Li{f_{k}}(-t)},e^{t}-1\right)$, $B_{n}^{(s)}(x)\sim\left(\left(\frac{e^{t}-1}{t}\right)^{s},t\right)$, let us assume that
\begin{equation}\label{eq:50}
\tilde{A}_{n}^{(r,k)}(x)=\sum_{m=0}^{n}C_{n,m}B_{m}^{(s)}(x),\,\,\,\,(r,s\in\mathbf{N}).
\end{equation}
By (\ref{eq:21}), we get
\begin{align}\label{eq:51}
C_{n,m}&=\frac{1}{m!}\Bigg\langle\left(\frac{t}{(1+t)\log{(1+t)}}\right)^{r}\left(\frac{t}{\log{(1+t)}}\right)^{s}Li{f_{k}}\left(-\log{(1+t)}\right)\\
&\quad\times\left(\log{(1+t)}\right)^{m}\Bigg\vert x^{n}\Bigg\rangle\nonumber\\
&=\frac{1}{m!}\sum_{l=0}^{n-m}\frac{m!}{(l+m)!}S_{1}(l+m,m)(n)_{l+m}\nonumber\\
&\quad\times\left\langle\left(\frac{t}{(1+t)\log{(1+t)}}\right)^{r+s}Li{f_{k}}\left(-\log{(1+t)}\right)(1+t)^{s}\Bigg\vert x^{n-l-m}\right\rangle\nonumber\\
&=\sum_{l=0}^{n-m}\binom{n}{l+m}S_{1}(l+m,m)\left\langle\sum_{a=0}^{\infty}\tilde{A}_{a}^{(r+s,k)}(s)\frac{t^{a}}{a!}\Bigg\vert x^{n-l-m}\right\rangle\nonumber\\
&=\sum_{l=0}^{n-m}\binom{n}{l}S_{1}(n-l,m)\tilde{A}_{l}^{(r+s,k)}(s).\nonumber
\end{align}
Therefore, by (\ref{eq:50}) and (\ref{eq:51}), we obtain the following theorem.

\begin{thm}\label{eq:thm7}
For $n\geq 0$, $r,s\in\mathbf{N}$, we have
\begin{equation*}
\tilde{A}_{n}^{(r,k)}(x)=\sum_{m=0}^{n}\left\{\sum_{l=0}^{n-m}\binom{n}{l}S_{1}(n-l,m)\tilde{A}_{l}^{(r+s,k)}(s)\right\}B_{m}^{(s)}(x).
\end{equation*}
\end{thm}

\noindent For $\tilde{A}_{n}^{(r,k)}(x)\sim\left(\left(\frac{te^{t}}{e^{t}-1}\right)^{r}\frac{1}{Li{f_{k}}(-t)},e^{t}-1\right)$, $H_{n}^{(s)}(x\vert\lambda)\sim\left(\left(\frac{e^{t}-\lambda}{1-\lambda}\right)^{s},t\right)$, let us assume that
\begin{equation}\label{eq:52}
\tilde{A}_{n}^{(r,k)}(x)=\sum_{m=0}^{n}C_{n,m}H_{m}^{(s)}\left(x\vert\lambda\right),
\end{equation}
where $\lambda\in\mathbf{C}$ with $\lambda\neq 1$, $r,s\in\mathbf{N}$ and $k\in\mathbf{Z}$.\\
From (\ref{eq:21}), we have
\begin{align}\label{eq:53}
C_{n,m}&=\frac{1}{m!}\Bigg\langle\left(\frac{t}{(1+t)\log{(1+t)}}\right)^{r}Li{f_{k}}\left(-\log{(1+t)}\right)\left(1+\frac{t}{1-\lambda}\right)^{s}\\
&\quad\times\left(\log{(1+t)}\right)^{m}\Bigg\vert x^{n}\Bigg\rangle\nonumber\\
&=\frac{1}{m!}\sum_{l=0}^{n-m}\frac{m!}{(l+m)!}S_{1}(l+m,m)(n)_{l+m}\nonumber\\
&\quad\times\left\langle\left(\frac{t}{(1+t)\log{(1+t)}}\right)^{r}Li{f_{k}}\left(-\log{(1+t)}\right)\left(1+\frac{t}{1-\lambda}\right)^{s}\Bigg\vert x^{n-l-m}\right\rangle\nonumber\\
&=\sum_{l=0}^{n-m}\binom{n}{l+m}S_{1}(l+m,m)\sum_{a=0}^{n-l-m}\binom{s}{a}\left(\frac{1}{1-\lambda}\right)^{a}(n-l-m)_{a}\nonumber\\
&\quad\times\left\langle\left(\frac{t}{(1+t)\log{(1+t)}}\right)^{r}Li{f_{k}}\left(-\log{(1+t)}\right)\Bigg\vert x^{n-l-m-a}\right\rangle\nonumber\\
&=\sum_{l=0}^{n-m}\binom{n}{l+m}S_{1}(l+m,m)\sum_{a=0}^{n-l-m}\binom{s}{a}\left(\frac{1}{1-\lambda}\right)^{a}(n-l-m)_{a}\tilde{A}_{n-l-m-a}^{(r,k)}\nonumber\\
&=\sum_{l=0}^{n-m}\sum_{a=0}^{l}\frac{\binom{n}{l}\binom{s}{a}\binom{l}{a}a!}{\left(1-\lambda\right)^{a}}S_{1}(n-l,m)\tilde{A}_{l-a}^{(r,k)}.\nonumber
\end{align}
Therefore, by (\ref{eq:52}) and (\ref{eq:53}), we obtain the following theorem.

\begin{thm}\label{eq:thm8}
For $n\geq 0$, $r,s\in\mathbf{N}$, $k\in\mathbf{Z}$ and $\lambda\in\mathbf{C}$ with $\lambda\neq 1$, we have
\begin{equation*}
\tilde{A}_{n}^{(r,k)}(x)=\sum_{m=0}^{n}\left\{\sum_{l=0}^{n-m}\sum_{a=0}^{l}\frac{\binom{n}{l}\binom{s}{a}\binom{l}{a}a!}{\left(1-\lambda\right)^{a}}S_{1}(n-l,m)\tilde{A}_{l-a}^{(r,k)}\right\}H_{m}^{(s)}\left(x\vert\lambda\right).
\end{equation*}
\end{thm}

\noindent Let us consider the following two Sheffer sequences:
\begin{equation*}
\tilde{A}_{n}^{(r,k)}(x)\sim\left(\left(\frac{te^{t}}{e^{t}-1}\right)^{r}\frac{1}{Li{f_{k}}(-t)},e^{t}-1\right)
\end{equation*}
and
\begin{equation*}
(x)_{n}=\sum_{l=0}^{n}S_{1}(n,l)x^{l}\sim\left(1,e^{t}-1\right).
\end{equation*}
Suppose that
\begin{equation}\label{eq:54}
\tilde{A}_{n}^{(r,k)}(x)=\sum_{m=0}^{n}C_{n,m}(x)_{m}.
\end{equation}
By (\ref{eq:21}), we get
\begin{equation}\label{eq:55}
C_{n,m}=\binom{n}{m}\tilde{A}_{n-m}^{(r,k)}.
\end{equation}
Therefore, by (\ref{eq:54}) and (\ref{eq:55}), we get
\begin{equation}\label{eq:56}
\tilde{A}_{n}^{(r,k)}(x)=\sum_{m=0}^{n}\binom{n}{m}\tilde{A}_{n-m}^{(r,k)}(x)_{m},
\end{equation}
where $r,n\geq 0$ and $k\in\mathbf{Z}$.


\bigskip
ACKNOWLEDGEMENTS. This work was supported by the National Research Foundation of Korea(NRF) grant funded by the Korea government(MOE)\\
(No.2012R1A1A2003786 ).
\bigskip


\noindent
\author{Department of Mathematics, Sogang University, Seoul 121-742, Republic of Korea
\\e-mail: dskim@sogang.ac.kr}\\
\\
\author{Department of Mathematics, Kwangwoon University, Seoul 139-701, Republic of Korea
\\e-mail: tkkim@kw.ac.kr}
\end{document}